\documentclass{amsart}

\usepackage{a4wide}
\usepackage{amsmath}
\usepackage{amsfonts}
\usepackage{rotating}
\usepackage{bbm}
\usepackage[all,cmtip]{xy}

\usepackage{amssymb}
\usepackage[utf8]{inputenc}
\usepackage{graphicx}
\usepackage{dsfont}
\usepackage{color}

\newtheorem{theorem}{Theorem}[section]

\newtheorem{proposition}[theorem]{Proposition}
\newtheorem{corollary}[theorem]{Corollary}

\theoremstyle{definition}
\newtheorem{definition}[theorem]{Definition}

\theoremstyle{remark}
\newtheorem{remark}[theorem]{Remark}

\numberwithin{equation}{section}

\newcommand{\R}{{\mathbb{R}}}

\newcommand{\Z}{{\mathbb{Z}}}
\newcommand{\N}{{\mathbb{N}}}

\renewcommand{\epsilon}{\varepsilon}
\renewcommand{\phi}{\varphi}
\renewcommand{\theta}{\vartheta}

\newcommand{\w}{\wedge}

\title{Expectation of topological invariants}
\author{Taejin Paik}
\address{Department of Mathematical Sciences and Research Institute of Mathematics, Seoul National University\\
Building 27, room 439\\
San 56-1, Sillim-dong, Gwanak-gu, Seoul, South Korea\\
Postal code 08826 }

\email{paiktj@snu.ac.kr}

\author{Otto van Koert}
\address{Department of Mathematical Sciences and Research Institute of Mathematics, Seoul National University\\
Building 27, room 402\\
San 56-1, Sillim-dong, Gwanak-gu, Seoul, South Korea\\
Postal code 08826 }

\email{okoert@snu.ac.kr}

\date{}

\newcommand{\mbb}[1]{\mathbb{#1}}
\newcommand{\mcal}[1]{\mathcal{#1}}

\begin{document}

\begin{abstract}
In this paper, we study the expectation values of topological invariants of the Vietoris-Rips complex and Čech complex for a finite set of sample points on a Riemannian manifold.
We show that the Betti number and Euler characteristic of the complexes are Lipschitz functions of the scale parameter and that there is an interval such that the Betti curve converges to the Betti number of the underlying manifold.
\end{abstract}

\maketitle

\section{Introduction}
An important result in topological data analysis concerns the stability theorem for persistent homology; one version states that the barcodes of the persistent homologies of sublevel sets of two continuous functions that are close must be close in bottleneck distance \cite{stability}.
Such a result is necessary for topological data analysis.
For example, if we associate the Vietoris-Rips complex with a point cloud, then the persistent homology of this complex will depend on the positions of the points in the point cloud. The stability theorem tells us that a small perturbation of these points will not change the barcodes in an essential way, thus giving us stable results for this type of noise.

On the other hand, the persistent homology of the Vietoris-Rips complex is not stable for outliers or the addition of data points. This has been observed in several papers. In \cite{BGMP} this phenomenon has been given a mathematical definition and referred to as fragility. 
A simple explicit example of this phenomenon is obtained by taking point samples near a circle, say a noisy circle.
Adding a single data point near the center will cause an enormous change in the barcodes. 
In short, persistent homology is not stable for this type of noise.
Although outliers are, by definition, uncommon, they do occur, so this phenomenon must be dealt with if one wishes to obtain meaningful statistical inference about the shape of the data. There are already various papers dealing with these issues, see \cite{BGMP, Bubenik}, but we will propose a different approach below.

Another issue concerns one of the motivating goals of topological data analysis; namely, we want to find the shape of data. 
The Vietoris-Rips complex achieves this in the following sense; a result of Latschev, \cite{VR_homotopy}, tells us that for sufficiently small scale parameter, the Vietoris-Rips complex of many evenly spaced out sample points on a Riemannian manifold $M$ reproduces the homotopy type of $M$.
Unfortunately, if the scale parameter is too large, the resulting Vietoris-Rips complex can have a completely different homotopy type.
For example, the Vietoris-Rips complex of the circle can be homotopy equivalent to any odd-dimensional sphere depending on the scale parameter, \cite{Ada1}.
This is, of course, highly undesirable for an unknown data set where we do not know what scale parameters to choose.
Although persistent homology deals with this problem of choice by capturing different scales, the persistent homology of the Vietoris-Rips complex will still have the above phenomenon.

In this paper, we investigate these issues in a simple setting; we consider the expectation value of a topological invariant of the Vietoris-Rips complex associated with a point cloud $\mbb{X} = \{ X_i \}$ lying on a manifold $M$.
Several invariants, such as the Betti number, can be used. For flexibility, we consider a topological invariant $T$ of a simplicial complex whose numerical value is bounded in terms of the number of simplices, such as the Betti numbers or Euler characteristic (see Section~\ref{sec:proof_continuity} for more detail).
After taking independent, identically distributed sample points $\mbb{X}=\{ X_i \}$ on $M$ using a probability distribution, we can view $T(\mbb{X}, t):= T(\operatorname{VR}(\mbb{X}, t))$ as a random variable.
We prove that the expectation value of this random variable is Lipschitz continuous in $t$ under some assumptions on the manifold $M$.
More precisely, we have the following theorem.
\begin{theorem}
\label{thm:continuity_expectation}
Suppose that $M$ is an oriented, complete Riemannian manifold equipped with a probability measure $\mu$ induced by a semi-positive top-form. 
We associate a random variable $T(\mbb{X},t)$ with $n$ independent, identically distributed samples $\mbb{X}=(X_1,\ldots,X_n)$ on $M$. 
Then the resulting function 
$$
t \longmapsto \mbb{E}[T(\mbb{X},t)]
$$
is continuous in $t$. If, in addition, there is a constant $K$ such that $|\mu( B_b(p) ) -\mu(B_a(p))|\leq K |b-a|$ for every $p$ in $M$ and $a,b$ in $\R$, then $t\mapsto \mbb{E}[T(\mbb{X},t)]$ is Lipschitz continuous.
\end{theorem}
The continuity of this expectation value means that the statistics of the invariant obtained by repeated sampling are stable for the scale parameter $t$ of the Vietoris-Rips complex.
In practice, the graph of expectation values displays even smooth behavior since Lipschitz functions are smooth almost everywhere by Rademacher's theorem \cite{lip_almost}.
The continuity of the expectation together with the next result shows that we have stability for outliers provided we measure sufficiently often, or perform a procedure such as bootstrapping. 
\begin{theorem}
\label{thm:moments}
Assume that $M$ is a compact, oriented Riemannian manifold with probability measure $\mu$ induced by a positive top-form.
Let $S(\mathbb{X}_n,t)$ denote either the Čech complex $\check{C}$ or the $\operatorname{VR}$-complex associated with $n$ independent, identically distributed samples $X_1,\ldots, X_n$ in $M$, and consider a topological invariant $T$ of a simplicial complex whose numerical value is bounded by a polynomial in the number of sample points $n$, such as the Betti numbers. We abbreviate
$$
T(\mathbb{X}_n,t):=T(S(\mathbb{X}_n,t) ).
$$
Then there is $\epsilon_0>0$ such that the following holds for $t\in (0,\epsilon_0)$
\begin{itemize}
\item The expectation $\mbb{E}[T(\mathbb{X}_n,t)]$ converges to $T(M)$ as $n$ goes to infinity.
\item The variance $\mbb{E}[T(\mathbb{X}_n,t)^2] - \mbb{E}[T(\mathbb{X}_n,t)]^2$ converges to $0$ as $n$ goes to infinity.
\end{itemize}
\end{theorem}
Similar results can also be found in the literature. We mention, for example, \cite{niyogi}.
Our approach is somewhat different and focuses on expectation values.
We hope our estimates can be refined if more information is known about the manifold of interest.
To support this hope, we explicitly compute the expectation values and variances of the first Betti number of $S^1$ with a uniform probability measure as a function of the scale parameter in Section~\ref{sec:betti_circle}.
A weakness in the above statement is that we don't have good a priori information on $\epsilon_0$.
As the proof explains, the scale $\epsilon_0$ depends on the convexity radius, which cannot be determined from point samples only, and this radius can be very small.
On the positive side, the theorem shows that there is an interval of small $t$-parameters for which the expectation value gets close to the true value of the topological invariant.  

\subsection*{Acknowledgements}
Taejin Paik was supported by National Research Foundation of Korea (NRF) Grant funded by the Korean Government, No. 2017R1A5A1015626.

\section{Preliminaries}

\subsection{Simplicial complexes and notation}
\label{sec:cech_vr_samples}
We will briefly review some of the various simplicial complexes that we can associate with metric spaces $(X,d)$.
The {\bf Čech complex} at scale $t$ is defined as the simplicial complex
$$
\check{C}(X,t)
= \{
\sigma = (x_0,\ldots,x_k) ~|~\cap_{j=0}^k B_t(x_j) \neq \emptyset \}
.
$$
Put differently, the Čech complex at scale $t$ is defined as the nerve of the cover $\{ U_x = B_t(x) \}_{x\in X}$.
The Vietoris-Rips complex is a related concept.
An efficient definition is the following. The {\bf Vietoris-Rips complex} at scale $t$ is the simplicial complex defined by
$$
\operatorname{VR}(X,t)=\{ \emptyset \neq \sigma  \subseteq  X ~|~\mathrm{diam}(\sigma) \leq t \}
.
$$
From the definitions, we  see that $\check{C}(X,r) \subseteq \operatorname{VR}(X,2r) $ and $\operatorname{VR}(X,2r)  \subseteq \check{C}(X,r+\epsilon) $ for all $\epsilon>0$.

We will consider a manifold $M$ of dimension at least $1$ equipped with a distance function $d$.
Given $n$ not necessary distinct points $X_1,\ldots,X_n$ on $M$, we get a finite metric space $(\mathbb{X}=\{ X_1,\ldots,X_n \},d|_{\mathbb{X}} )$ where $d|_{\mathbb{X}}$ is the (semi)-distance function induced on $\mathbb{X}$.
With some abuse of notation, we will also write $\mathbb{X}=(X_1,\ldots,X_n)$ for the $n$-tuple of points in $M^n=M\times \ldots \times M$.

In the next section, we will review some notions regarding probability measures on manifolds.
After introducing a probability measure on $M^n$, we will treat the $n$-tuple $\mathbb{X}$ as a random variable, and we will investigate the expectation value of topological invariants such as the Betti numbers of the resulting random $\operatorname{VR}$-complex as a function of the scale parameter.
For this purpose, it is helpful to list our notation to keep track of the parameter values where these simplicial complexes change.
For readability, we also give some names to the sets.
\begin{itemize}
\item We denote the number of edges in $\operatorname{VR}(\mbb{X},t)$ by $e(\mbb{X},t)$.
\item For $k=0,1, \ldots, {n \choose 2}$ we define 
$$
\ell_k(\mbb{X}):=\min\{ t\in \R_{\geq 0} ~|~e(\mbb{X},t) \geq k \}.
$$
In words, this is the smallest scale parameter $t$ for which $\operatorname{VR}(\mbb{X},t)$ has $k$ edges.
We will say that an edge has been created at scale $\ell_k(\mbb{X})$ for the first time. From this definition, we immediately obtain the following
$$
\ell_{e(\mathbb{X},t)}(\mbb{X}) \leq t < \ell_{e(\mathbb{X},t)+1}(\mbb{X}).
$$
\item We define the edge creation set 
$$
eC_n^M(t):=
\{
\mathbb{X}=(X_1,\ldots, X_n) \in M^n~|~\ell_{e(\mathbb{X},t)}(\mathbb{X})=t
\}.
$$
This set describes configurations of points for which an edge has been created at scale $t$ for the first time.
\item We define the set of pairs as
$$
P_n^M:=\{ \mathbb{X}=(X_1,\ldots,X_n)\in M^n ~|~\ell_k(\mathbb{X}) = \ell_{k+1}(\mathbb{X}) \text{ for some }k \in \Z_{>0} \}.
$$
This set describes configurations for which there are at least two pairs of points among $\mathbb{X}$ with equal distance.
\item We define the set of collisions as
$$
Col_n^M= \{ (X_1,\ldots,X_n)\in M^n~|~ X_i=X_j \text{ for some }i\neq j \}.
$$
\item With a given $n$-tuple of points $\mbb{X}=(X_1,\ldots,X_n)$, we associate its {\bf covering radius} as
$$
\rho_M(\mbb{X}):=
\inf
\left\{
r \in \R_{\geq 0}~\bigg|~ \bigcup_{i=1}^n B_r(X_i) = M
\right\}
.
$$
\end{itemize}

We will use the following two theorems to connect the Čech complex and Vietoris-Rips complex back to the homotopy type of $M$.
The nerve theorem is a classical result and can, for instance, be found in \cite[Chapter 4.G]{hatcher}. For the Vietoris-Rips complex, we have a result due to Latschev, \cite{VR_homotopy}.
\begin{theorem}[Nerve theorem]
\label{thm:nerve}
Suppose that $M$ is a paracompact space with an open cover $\mcal{U} = \{U_\alpha\}$. If every nonempty intersection of finitely many $U_\alpha$ is contractible, then $M$ is homotopy equivalent to the nerve of $\mcal{U}$.
\end{theorem}

\begin{theorem}[Homotopy type of Vietoris-Rips complex]
\label{thm:homotopy_vr}
Let $M$ be a closed Riemannian manifold.
Then there exists $\epsilon_0 > 0$ such that for every $0<\epsilon \leq \epsilon_0$ there exists a $\delta > 0$ such that if Gromov-Hausdorff distance between $M$ and a finite subset $\mbb{X}$ of $M$ is less than $\delta$ then $\operatorname{VR}(\mbb{X}, \epsilon)$ is homotopy equivalent to $M$.
\end{theorem}

\subsection{Probability measures on manifolds}
\label{sec:prob_mfd}
We will collect some facts from the literature.
Consider a smooth, oriented manifold $M$ equipped with a continuous top-form $\omega$.
We call $\omega$ {\bf semi-positive} if for every oriented coordinate chart $\phi: V \subset \R^n \to U \subset M$ we have
$$
\phi^* \omega = gdx^1 \w \ldots \w dx^n
$$
with $g\geq 0$. If this inequality is strict, i.e.~if $g>0$, then we call $\omega$ a {\bf positive} top-form.
As explained in \cite[Chapter XVI, paragraph 4]{Lang:DG}, for example, a positive top-form gives rise to a measure $\mu$. A similar construction works for semi-positive top-forms.

In the computations below, we can exclude the contribution of certain subsets to the integrals we need to compute by measure zero arguments. These conclusions are usually drawn using the following version of measure zero commonly used in differential topology.
\begin{definition}
Let $M$ be a smooth $n$-dimensional manifold, possibly with boundary.
We say that a subset $A\subseteq M$ has {\bf measure zero in geometric sense} if for every smooth chart $(U,\phi)$ in $M$, the subset $\phi(A \cap U) \subseteq \R^n$ has measure zero in Lebesgue sense. 
\end{definition}
On the other hand, we could also directly define measure zero on a manifold using the above measure $\mu$. Under mild assumptions, these two concepts coincide. We have
\begin{theorem}
Suppose that $M$ is an oriented manifold equipped with a measure $\mu$ induced by a semi-positive top-form $\omega$. If a subset $A$ in $M$ has measure zero in geometric sense, then $\mu(A)=0$. 
For the converse, assume in addition that $\omega$ is a positive top-form. If with this additional assumption, we have $\mu(A) =0$ for some subset $A$ in $M$, then $A$ has measure zero in geometric sense.
\end{theorem}
We omit the proof. Below we will assume that we have a probability measure $\mu$ on a manifold $M$ induced by a semi-positive top-form $\omega$, so $\mu(M)=1$. This gives us an induced probability measure $\mu^n$ on $M^n=M\times \ldots \times M$. 
We view an $n$-tuple $(X_1,\ldots,X_n)\in M^n$ as $n$ independent, identically distributed samples taken on $M$.

\subsection{Ingredients from Riemannian geometry}
We collect some ingredients from Riemannian geometry. Good references are \cite{CE, Sakai}.
Suppose that $(M,g)$ is a Riemannian manifold. We obtain a distance function by the formula
$$
% d(p,q)= \inf_{\gamma \in C^\infty([0,1],M), \gamma(0)=p, \gamma(1)=q} len(\gamma),
d(p,q)= \inf_{\substack{\gamma \in C^\infty([0,1],M),\\ \gamma(0)=p, \gamma(1)=q}} len(\gamma),
\text{ where }
len(\gamma)=\int_{0}^1 g_{\gamma(t)}(\dot \gamma(t),\dot \gamma(t))dt.
$$
As a shortcut notation, we write $d_p(x):=d(p,x)$ whenever the point $p$ is fixed.
We also recall the exponential map. Given $p\in M$, this is the map
$$
\exp_p:T_p M \longrightarrow M,\quad v \longmapsto \gamma_{p,v}(1),
$$
where $\gamma_{p,v}$ is the geodesic with $\gamma_{p,v}(0)=p$ and $\dot\gamma_{p,v}(0)=v$. 
Now fix a unit vector $v\in T_pM$.
By formula~(4.1) in chapter III of \cite{Sakai}, we know that 
$$
d(p,\exp_p(tv))=t
$$
for small $t$.
If $t_0=\sup\{ t~|~d(p,\exp_p(tv))=t \}$ is finite, then we call $\gamma_{p,v}(t_0)$ a cut point of $\gamma_{p,v}$ with respect to $p$.
If $M$ is compact, then every geodesic ray has a cut point.
The {\bf cut locus} $C(p)$ of a point $p$ is the set of cut points with respect to $p$.
We collect some facts about the cut locus. From the book of Sakai, \cite{Sakai}, Chapter III, Lemma~4.4 and Proposition~4.8, we have.
\begin{theorem}
The cut locus $C(p)$ of a point $p$ in a Riemannian manifold $(M,g)$ has measure zero. 
\end{theorem}
\begin{proposition}
Suppose that $p$ is a point in a complete Riemannian manifold $(M,g)$. 
Then the distance function $d$ is smooth on $M \setminus ( C(p) \cup \{ p \} )$. 
\end{proposition}
We will denote the restriction $d_p|_{M \setminus ( C(p) \cup \{ p \} )}$ by $\hat d_p$.
We observe that the function $\hat d_p$ has no critical points.
Indeed, any $q$ can be written as $q=\exp_p(tv)$ where $v$ is a unit vector. By the above $\hat d_p(q)=\hat d_p(\exp_p(tv))=t$.
As the derivative with respect to $t$ is non-zero, we see that $\hat d_p$ has no critical points.
As a corollary, we obtain the following.
\begin{proposition}
\label{prop:dist_inv_measure}
Suppose that $p$ is a point in a complete, oriented Riemannian manifold $(M,g)$ equipped with a measure induced by a semi-positive top-form. Then for every $t$ in $\R$ the set $d_p^{-1}(t)$ has measure zero.
\end{proposition}
To see this, we note that 
$$
d_p^{-1}(t) \subset \hat d_p^{-1}(t) \cup C(p) \cup \{ p \} .
$$
As every value of $\hat d_p$ is a regular value, $\hat d_p^{-1}(t)$ is a submanifold of codimension $1$, so it has measure $0$. Together with the previous theorem, the claim follows.

\begin{definition}
We call a set $C$ in a Riemannian manifold $(M,g)$ {\bf strongly convex} if for every pair of points $p,q$ in the closure $\overline{C}$, there is unique minimal geodesic $\gamma$ connecting $p$ to $q$ whose interior is contained in $C$. 
The {\bf convexity radius} of a point $p$ in $M$ is then defined as
$$
cvxrad(p):=\sup \{ r > 0 ~|~B_s(p) \text{ is strongly convex for } 0<s<r \}
,
$$
The convexity radius of $(M,g)$ is then defined as $cvxrad_M:=\inf cvxrad(p)$.
\end{definition}

Combining this notion with the nerve theorem, we have the following direct corollary.
\begin{corollary}
\label{cor:convex_cech_he}
Suppose that $\mbb{X}=(X_1,\ldots,X_n)$ is an $n$-tuple of points on a Riemannian manifold $M$.
Assume that we are given a real number $\rho_0$ such that $\rho_M(\mbb{X}) <\rho_0<cvxrad_M$.
Then the Čech complex $\check{C}(\mbb{X},\rho_0)$ is homotopy equivalent to $M$.
\end{corollary}

\section{Proof of the continuity theorem}
\label{sec:proof_continuity}
We start by showing that the set of parameters and configurations where the $\operatorname{VR}$-complex changes has measure $0$.
More precisely, we have
\begin{proposition}
\label{prop:measure_level_set_d_0}
Suppose that $(M,g)$ is an oriented, complete Riemannian manifold. Assume that $\mu$ is a probability measure on $M$ induced by a semi-positive top-form. As in Section~\ref{sec:prob_mfd}, we denote the induced measure on $M^n$ by $\mu^n$.
Then the edge creation set $eC^M_n(t)$, the set of pairs $P_n^M$ and the set of collisions $Col_n^M$ all have measure $0$,
$$
\mu^n(eC^M_n(t) )= \mu^n(P_n^M )= \mu^n(Col_n^M) =0.
$$
\end{proposition}
\begin{proof}
We begin with the set of collisions. It suffices to show that $\mu^n(X_1=X_2)=0$.
By Fubini's theorem and Proposition~\ref{prop:dist_inv_measure} we obtain
\[
\begin{split}
\mu^n(\{\mathbb{X}~|~X_1=X_2\}) & = \mu^2(\{ (X_1,X_2)\in M^2~|~X_1=X_2 \} ) \\
&= \int_M \int_M \chi_{d_{X_2}^{-1}(0)} d\mu_{X_1} d\mu_{X_2} \\
&=\int_{M} \mu(d_{X_2}^{-1}(0)) d\mu_{X_2} =0.
\end{split}
\]
Now we will show that the set of pairs has measure $0$. We need to show the following two cases.
\[
\begin{split}
&\mu^n( \{ \mathbb{X}~|~0=\ell_0(\mathbb{X})=\ell_1(\mathbb{X}) \} ) = 0, \text{ and }\\
&\mu^n( \{ \mathbb{X}~|~\ell_k(\mathbb{X})=\ell_{k+1}(\mathbb{X}) \text{ for some }k\neq 0 \} ) = 0.
\end{split}
\]
For the first case, note that if $\ell_1$ vanishes, then there are $i, j$ with $i\neq j$ such that $X_i=X_j$.
As we already know that $\mu^n(Col_n^M)=0$, the first case holds.
In the second case, we have $\ell_k=\ell_{k+1}$, so there are two pairs of points with the same distance. This occurs in the following situations,
$$
d(X_{i_1},X_{i_2})=d(X_{i_3},X_{i_4}) \text{ or }d(X_{i_1},X_{i_2})=d(X_{i_2},X_{i_3}) 
$$
for indices $i_1,i_2,i_3,i_4$ that are all distinct. 
We compute the measure of the set corresponding to the first situation.
This equals
\[
\begin{split}
\mu^n(
\{ 
\mathbb{X}~|~d(X_{i_1},X_{i_2})=d(X_{i_3},X_{i_4})
\}
)
&=\mu^n(
\{ 
\mathbb{X}~|~d_{X_2}(X_1)=d_{X_3}(X_{4})
\}
)\\
&=\mu^4(
\{ 
\mathbb{X}~|~d_{X_2}(X_1)=d_{X_3}(X_{4})
\}
) \\
&= \int_{M^4} \chi_{d_{X_2}^{-1}(d_{X_3}(X_4))} d\mu_{X_1} d\mu_{X_2} d\mu_{X_3} d\mu_{X_4}.
\end{split}
\]
As in the previous case, this vanishes because of Fubini's theorem and Proposition~\ref{prop:dist_inv_measure}. The second situation, i.e.~when $d(X_{i_1},X_{i_2})=d(X_{i_2},X_{i_3})$ can be dealt with in the same way.
Finally, we address the edge creation set $eC_n^M(t)$ for $t>0$.
By definition $\mathbb{X}$ is in $eC_n^M(t)$ precisely when we have $X_i,X_j$ such that $d(X_i, X_j)=t$.
In order to show that $\mu^n( eC_n^M(t) )$ vanishes, it therefore suffices to show that 
$$
\mu^2(
\{
(X_1,X_2) \in M^2 ~|~d(X_1,X_2)=t  
\}
)=0.
$$
Again, this follows in the same way as before from Fubini's theorem and Proposition~\ref{prop:dist_inv_measure}.
\end{proof}

Given an $n$-tuple of points $\mathbb{X}$ on a manifold $M$, we construct the $\operatorname{VR}$-complex at scale $t$. 
We consider a topological invariant $T$ of a simplicial complex satisfying the following:
there is $f: \N \to \R$ such that $|T(S)| \leq f( |S| )$ for every simplicial complex $S$.
The Euler characteristic and the $i$-th Betti number satisfy this condition.
Theorem~\ref{thm:continuity_expectation}, which we prove below, tells us the function $t\longmapsto \mbb{E}[T(\operatorname{VR}(\mathbb{X},t))]$ is Lipschitz continuous.
To make the notation less bloated, we will abbreviate $T(\mathbb{X},t):=T(\operatorname{VR}(\mathbb{X},t))$.

\begin{proof}[Proof of Theorem~\ref{thm:continuity_expectation}]
We start with some observations.
By our assumption on $T$, there is $f: \N \to \R$ such that $|T(S)| \leq f( |S| )$ for every simplicial complex $S$.
The number of simplices in the $\operatorname{VR}$-complex of $n$ points is bounded from above by $\sum_{k=0}^n {n \choose k}=2^n$, so $T(\mathbb{X},t) \leq f(2^n)$.
For $0\leq a \leq b$, we define the set
$$
S(a,b):= \{ 
\mbb{X}~|~\text{there is }k \text{ such that } a<\ell_k(\mbb{X}) < b \}.
$$
This set describes configurations $\mbb{X}\in M^n$ such that the $\operatorname{VR}$-complexes $\operatorname{VR}(\mbb{X},a)$ and $\operatorname{VR}(\mbb{X},b)$ differ.
As $S(a,b)$ is open, it is measurable. Furthermore, we have the following property:
$$
[a_0,b_0] \subseteq [a_1,b_1] \text{ implies } S(a_0,b_0) \subseteq S(a_1,b_1).
$$
We now take $0\leq a < b$ and estimate the difference in expectation values using the above proposition and the definition of $S(a,b)$,
\[
\begin{split}
|\mbb{E}[T(\mbb{X},a)]-\mbb{E}[T(\mbb{X},b)]| & \leq \mbb{E}[|T(\mbb{X},a)- T(\mbb{X},b)|] \\
& = \int_{M^n}|T(\mbb{X},a)-T(\mbb{X},b)| d \mu^n \\
& =\int_{S(a,b)} |T(\mbb{X},a)-T(\mbb{X},b)| d \mu^n \\
&\leq 2 f(2^n) \mu^n(S(a,b)).
\end{split}
\] 
We will now bound the probability of $S(a,b)$. 
From the definition of $S(a,b)$, we see that a configuration $\mbb{X}$ lies in $S(a,b)$ precisely when there is an edge created at a scale between $a$ and $b$.
Hence there must be points $X_i,X_j$ with $a<d(X_i,X_j)<b$.
With this in mind, we compute
\[
\begin{split}
\mu^n(S(a,b)) &= \mu^n
\left(
\bigcup_{0<i<j\leq n} 
\{
\mbb{X}~|~a< d(X_i,X_j) <b 
\}
\right) \\
& \leq \sum_{0<i<j\leq n}\mu^n( \{ \mbb{X}~|~a<d(X_i,X_j)<b  \} ) \\
& =\sum_{0<i<j\leq n}\mu^2( \{ \mbb{X}~|~a<d(X_1,X_2)<b  \} ) \\
&= {n \choose 2} \int_M \mu( B_b(X_1) \setminus B_a(X_1))d\mu_{X_1}.
\end{split}
\]
We note that $\int_M \mu(B_b(X)) d\mu_X \leq \int_M d\mu_X =1$ and that $\mu( B_b(X) \setminus B_a(X))\leq \mu(B_b(X))$, which is an integrable function.
By the dominated convergence theorem, we see hence that
$$
\lim_{b \to a} \int_M \mu( B_b(X_1) \setminus B_a(X_1))d\mu_{X_1} =0,
$$ 
which shows that $\mbb{E}[T(\mbb{X},t)]$ is a continuous function of $t$.
For the final claim, we make that additional assumption that $|\mu( B_b(p) ) -\mu(B_a(p))|\leq K |b-a|$, which allows us to combine the above estimates to
$$
|\mbb{E}[T(\mbb{X},a)]-\mbb{E}[T(\mbb{X},b)]| \leq 2 \cdot {n \choose 2} \cdot f(2^n) K |b-a|.
$$
This means that $\mbb{E}[T(\mbb{X},t)]$ is Lipschitz.
\end{proof}

\begin{remark}
On compact manifolds, the additional assumption in the second part of the theorem,
$$
|\mu( B_b(p) ) -\mu(B_a(p))|\leq K |b-a|
$$
holds.
We also note that for two point samples on $S^1=\R/\Z$, the expectation value of the Euler characteristic is Lipschitz, but not smooth in $t=1/2$. Indeed, we have 
\[
\mbb{E}[\chi(\mathbb{X},t)]=
\begin{cases}
2(1-t) & \text{if }t\in [0,1/2] \\
1 & \text{if }t>1/2.
\end{cases}
\]
\end{remark}

\section{Expectation values and variance computations}
We consider a compact manifold $M$ with probability measure $\mu$ coming from a positive top-form.
Let $S(\mathbb{X}_n,t)$ denote either the Čech complex $\check{C}$ or the $\operatorname{VR}$-complex associated with the independent, identically distributed samples $\mathbb{X}_n$ in $M^n$.
We saw in the previous section that the expectation value of $T(\mathbb{X},t)$ is a Lipschitz function of the scale parameter $t$. 
Now we will impose the additional assumption that the topological invariant $T$ is polynomial in the number $n$ of sample points, and we will show that, for sufficiently small scale parameters, the expectation value of such a topological invariant $T$ converges to $T(M)$.
The main intuition for this comes from Corollary~\ref{cor:convex_cech_he}.

Given a scale parameter $\epsilon$, we first need to control the probability of the set of configurations of $n$ sample points $\mbb{X}$ on $M$ that do not yield a covering at scale $\epsilon$. 
\begin{proposition}
\label{prop:bad_set_small_measure}
Suppose that $M$ is a compact, oriented Riemannian manifold with a probability measure $\mu$ induced by a positive top-form, and take $\epsilon>0$.
Then the function defined by sending $n\mapsto\mu^n( \mbb{X}\in M^n~|~\rho_M(\mbb{X}) > \epsilon )$
is monotonically decreasing and converges to $0$ as $n$ goes to infinity.
\end{proposition}

\begin{proof}
Given an $n$-tuple of points $\mathbb{X}_n=(X_1,\ldots,X_n)$, we define for $k=1,\ldots,n$ the subtuple $\mathbb{X}_{k,n}:=(X_1,\ldots,X_k)$.
Since $\mu$ is a probability measure, we have
\[
\mu^{n}( \{ \mathbb{X}_{n}~|~ \rho_M(\mathbb{X}_{n}) >\epsilon \} )
=
\mu^{n+1}( \{ \mathbb{X}_{n+1}~|~ \rho_M(\mathbb{X}_{n,n+1}) >\epsilon \} )
\] 
Obviously, if the subset $\mathbb{X}_{n,n+1}$ suffices for an $\epsilon$-cover, then so does the full set $\mathbb{X}_{n+1}$.
In other words, if $\rho_M(\mathbb{X}_{n,n+1}) \leq \epsilon$, then $\rho_M(\mathbb{X}_{n+1})\leq \epsilon$.
This implies that
$$
\{ 
\mathbb{X}_{n+1}~|~\rho_M(\mathbb{X}_{n+1} ) >\epsilon
\}
\subseteq
\{ 
\mathbb{X}_{n+1}~|~\rho_M(\mathbb{X}_{n, n+1} ) >\epsilon
\}
,
$$
so taking their measure, we see that
\[
\begin{split}
\mu^{n+1}(\{ 
\mathbb{X}_{n+1}~|~\rho_M(\mathbb{X}_{n+1} ) >\epsilon
\})
&\leq
\mu^{n+1}(\{ 
\mathbb{X}_{n+1}~|~\rho_M(\mathbb{X}_{n, n+1} ) >\epsilon
\})\\
&=\mu^n(\{ 
\mathbb{X}_{n}~|~\rho_M(\mathbb{X}_{n} ) >\epsilon
\}),
\end{split}
\]
which means that $\mu^n(\{ 
\mathbb{X}_{n}~|~\rho_M(\mathbb{X}_{n} ) >\epsilon
\})$ is a monotonically decreasing function of $n$.

For the second statement, we choose a finite cover by $\epsilon/2$ balls, $M= \bigcup_{i=1}^k B_{\epsilon/2}(x_i)$, which we can do by the compactness of $M$.
We abbreviate $B_i:=B_{\epsilon/2}(x_i)$. 
We claim that for every $n$-tuple of sample points $\mathbb{X}_n$ with $\rho_M(\mathbb{X}_n)> \epsilon$, there is a ball $B_j$ that does not contain any of the $X_i$ for $i=1,\ldots,n$.
To see this, note that there is a point $x$ in $M \setminus \cup_{i=1}^nB_\epsilon(X_i)$ since the covering radius is larger. We find a ball $B_j$ containing this point, and by the triangle inequality, we see that 
$$
d(x_j,X_i) \geq d(x,X_i)-d(x,x_j) \geq \frac{\epsilon}{2}.
$$

With this claim in hand, we bound the probability that the covering radius is larger than $\epsilon$,
\[
\begin{split}
\mu^n(
\{ 
\mathbb{X}_{n}~|~\rho_M(\mathbb{X}_{n} ) >\epsilon
\}
)
&\leq 
\mu^n
(
\{ \mathbb{X}_n ~|~\text{there is }j\text{ such that } B_j \cap \mathbb{X}_n =\emptyset \}
) \\
&
\leq
\sum_{j=1}^k
\mu^n
(
\{
\mathbb{X}_n~|~ B_j \cap \mathbb{X}_n =\emptyset 
\}
) \\
&\leq \sum_{j=1}^k (1-\mu(B_j))^n \\
& \leq k(1-\mu(B_J) )^n,
\end{split}
\]
where we define $B_J$ as a minimizer of $\{ \mu(B_j)~|~j=1,\ldots,k \}$. Because $\mu$ is assumed to be induced by a positive top-form, $\mu(B_J)>0$, so the second claim follows.
\end{proof}

\begin{proof}[Proof of Theorem~\ref{thm:moments}]
We will show the stronger claim that all moments $T^i(\mathbb{X}_n,t)$ converge to $T^i(M)$ as $n$ goes to infinity provided $t$ is sufficiently small.
The arguments for the Čech and Vietoris-Rips complex are very similar and only differ in our initial choices. 
For the Čech complex, choose $\epsilon_0:=cvxrad_M$, and set $\delta_t:=t$.
For the Vietoris-Rips complex, we apply Theorem~\ref{thm:homotopy_vr} to obtain these numbers.
Indeed, we find that there is $\epsilon_0>0$ such that for $t \in(0,\epsilon_0)$ there is $\delta_t$ with the property that if the covering radius $\rho_M(\mathbb{X}_n)$ is less than $\delta_t$, then $S(\mathbb{X}_n,t)$ is homotopy equivalent to $M$.

Choose $t\in (0,\epsilon_0)$. Splitting the domain of integration $M^n$, we see that
\[
\mbb{E}[T^i(\mathbb{X}_n,t)]
=
\int
_{\rho_M(\mathbb{X}_n) <\delta_t} T^i(\mathbb{X}_n,t) d\mu^n
+
\int
_{\rho_M(\mathbb{X}_n) \geq \delta_t} T^i(\mathbb{X}_n,t) d\mu^n
\]
By the choice of $\epsilon_0$, and the condition $\rho_M(\mathbb{X}_n) <\delta_t$ we know that $S(\mathbb{X}_n,t)$ is homotopy equivalent to $M$ so the first term equals
\[
\int
_{\rho_M(\mathbb{X}_n) < \delta_t} T^i(\mathbb{X}_n,t) d\mu^n
=T^i(M) \cdot \mu^n( \{ 
\mathbb{X}_n ~|~\rho_M(\mathbb{X}_n) <\delta_t
\}).
\]
By assumption $T$ is polynomial in $n$, so there is some $j \in \N$ such that $T(\mathbb{X}_n,t)\leq C n^j$.
We use this to bound the second term in the above expectation value.
We find
\[
\begin{split}
\left\vert
\int_{\rho_M(\mathbb{X}_n)\geq \delta_t} T^i(\mathbb{X}_n,t) d\mu^n
\right\vert 
&
\leq
\int_{\rho_M(\mathbb{X}_n)\geq \delta_t} |T^i(\mathbb{X}_n,t)| d\mu^n \\
&\leq \mu^n(\rho_M(\mathbb{X}_n)\geq \delta_t) \cdot C^i n^{ij}.
\end{split}
\]
From Proposition~\ref{prop:bad_set_small_measure} and its proof we know that $
\mu^n(
\{ \mathbb{X}_n ~|~\rho_M(\mathbb{X}_n)\geq \delta_t \}
)
\leq k(1-\mu(B_J) )^n$,
where $k$ and $\mu(B_J)$ are positive real numbers depending only on $t$. Putting things together, we find
$$
\left\vert
\int_{\rho_M(\mathbb{X}_n)\geq \delta_t} T^i(\mathbb{X}_n,t) d\mu^n
\right\vert  \leq k(1-\mu(B_J) )^n C^i n^{ij}.$$ 
As $|1-\mu(B_J) |<1$, this converges to $0$ as $n$ goes to infinity.
\end{proof}

\section{Probability on a circle}
\label{sec:betti_circle}
If $\mbb{X} = \{X_i\}$ is a set of $n$ independent and identically distributed random samples of the uniform distribution on $(S^1=\R /\Z,g=dt\otimes dt)$ and $0<r<\frac{1}{3}$, then the probability that $\operatorname{VR}(\mbb{X}, r)$ is homotopy equivalent to $S^1$ is
\begin{equation}
\label{eq:circle_prob}
\mbb{P}(n, r) = nr^{n-1}\left[
\int_0^r g\left(n-1, \frac{1}{r}(1-x)\right) \; dx - rg\left(n-1, \frac{1}{r} - 1\right)
\right]
,
\end{equation}
where
\[
g(n, x) := \sum_{k=0}^{n}(-1)^{k}\binom{n}{k} (x-k)_{+}^{n}, \text{ with }
x_{+}= \begin{cases}x & x \geq 0 \\ 0 & x<0.\end{cases}
\]
We omit the computation of this probability, and only mention the key ingredients:
\begin{itemize}
\item a characterization when the Vietoris-Rips complex of finitely many points on $S^1$ is homotopy equivalent to $S^1$, which, together with stronger results, can be found in Corollary~4.5 of \cite{Ada1}.
\item the Irwin-Hall distribution, see Equation~(26.48) from \cite{JBK2}, describes the sum of independent, uniform random variables on $[0,1]$.
\end{itemize}

{

Using Equation~\eqref{eq:circle_prob} we obtain a formula for the expectation value and variance of the first Betti number. The results are plotted in Figure~\ref{fig:expectation_variance_plot}.
{
\begin{figure}[ht]
    \begin{center}
        \includegraphics[width=0.99\textwidth]{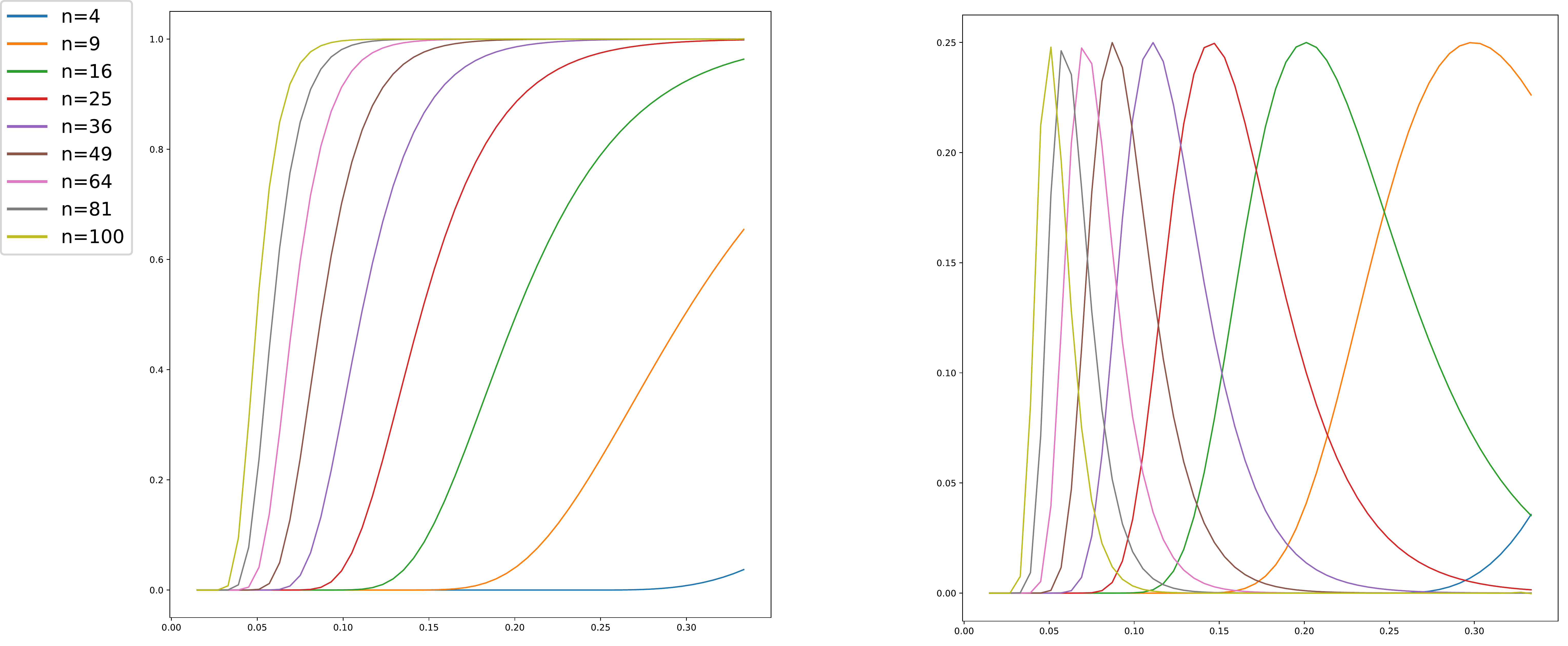}
        \caption{Expectation value and variance of the first Betti number of the VR complex of $n$ random points on the circle as obtained in Equation~\eqref{eq:circle_prob} }
        \label{fig:expectation_variance_plot}
    \end{center}  
\end{figure}
}

\end{document}